\documentclass[12pt]{article}
\usepackage{amsmath,amsthm,eucal}
\newcommand{\cl}{C\kern -0.2em \ell}
\newtheorem{theorem}{Theorem}
\newtheorem{prop}{Proposition}
\newcommand{\e}{\mbox{\bf e}}
\newcommand{\cA}{\mathcal{A}}
\newcommand{\p}{\prime}
\newcommand{\C}{\mbox{\bf C}} 
\newcommand{\R}{\mbox{\bf R}}
\newcommand{\re}{\mbox{\rm Re}\,}
\newcommand{\im}{\mbox{\rm Im}\,}
\newcommand{\M}{\mbox{\rm M}}
\newcommand{\BH}{{\rm\hskip 0.1pt %
                             I\hskip -2.15pt H}}
\begin{document}
\title{Dirac-Hestenes spinors and Weierstrass representation for surfaces
in 4D complex space}
\author{Vadim V. Varlamov}
\date{}
\maketitle
\begin{abstract}
Representations of Dirac-Hestenes and Dirac spinor fields via coordinates
of surfaces conformally immersed into 4-dimensional complex space are
proposed. A relation between time evolution of spinor fields and
integrable deformations of surfaces is discussed.
\end{abstract}
\section{Introduction}
For the first time Weierstrass representation for conformal immersion
of surface into $\R^3$ appeared in the result of variational problem on
search of minimal surface restricted by the some curve \cite{1}.
Generalization of Weierstrass formulae for surfaces with mean curvature
$H\neq 0$ was proposed by Eisenhart in 1909 \cite{2}. At present, the 
interest in this topic is increase after works of Kenmotsu \cite{3} and
Konopelchenko \cite{4}. An important relation between integrable
deformations of surfaces conformally immersed into $\R^3$ and nonlinear
differential equations of soliton theory has been established in paper \cite{4}.
Integrable deformations of surfaces (surfaces of revolution, surfaces of
constant mean curvature and so on), which are defined by hierarchies of
equations of nonlinear physics, considered in \cite{5}-\cite{10}.
A further generalization of Weierstrass representation onto a case of
4-dimensional real spaces with different signatures has been proposed by
Konopelchenko and Landolfi in recent papers \cite{11}-\cite{13}.

In present paper we consider Weierstrass representation for conformal
immersion of surfaces into 4-dimensional complex space. In section 2
generalized Weierstrass formulae in $\C^4$ are rewritten in a spinor
representation type form which coincides with the matrix representation of
a biquaternion algebra $\C_2\cong\M_2(\C)$ known in physics as a Pauli algebra.
It is allows using the well-known relation between Dirac-Hestenes and
Dirac spinors \cite{14,15} to establish a relation between coordinates for
surfaces immersed into $\C^4$ and Dirac spinors. An equivalence between
conjugated spinors and minimal right ideals is given in section 3.
Integrable deformations of surfaces in $\C^4$ defined by Davey-Stewartson
hierarchy and their relations with the time evolution of Dirac field
are considered in section 4.
\section{A spinor type form of Weierstrass representation for surfaces in
space $\C^4$}
Let $\C^4$ be a 4-dimensional complex space associated with a Dirac algebra
$\C_4$. Let us consider a generalized Weierstrass representation for
immersion of 2-dimensional surfaces in space $\C^4$. We propose that
generalized Weierstrass formulae in this case have a form \cite{16}
\begin{eqnarray}
X^1=\frac{i}{2}\int_\Gamma(\psi_1\psi_2d\bar{z}-\varphi_1\varphi_2dz),
\nonumber\\
X^2=\frac{1}{2}\int_\Gamma(\psi_1\psi_2d\bar{z}+\varphi_1\varphi_2dz),
\nonumber\\
X^3=\frac{1}{2}\int_\Gamma(\psi_1\varphi_2d\bar{z}-\varphi_1\psi_2dz),
\nonumber\\
X^4=\frac{i}{2}\int_\Gamma(\psi_1\varphi_2d\bar{z}+\varphi_1\psi_2dz),\label{e1}
\end{eqnarray}
where
\begin{equation}\label{e2}
\begin{array}{ccc}
\psi_\alpha&=&p\varphi_\alpha,\\
\varphi_\alpha&=&q\psi_\alpha,
\end{array}\quad\alpha=1,2.
\end{equation}
Here $\psi_\alpha,\,\varphi_\alpha$ and $p,\,q$ are complex-valued
functions, $\Gamma$ is a contour in complex plane $\C$. We will interpret
the functions $X^i(z,\bar{z})$ as the coordinates in $\C^4$. It is easy to
verify that components of induced metric have a form
\begin{eqnarray}
g_{zz}&=&\overline{g_{\bar{z}\bar{z}}}=\sum^4_{i=1}(X^i_z)^2=0,\nonumber\\
g_{z\bar{z}}&=&\sum^4_{i=1}(X^i_zX^i_{\bar{z}})=\psi_1\psi_2\varphi_1\varphi_2.
\nonumber
\end{eqnarray}
Therefore, the formulae (\ref{e1}), (\ref{e2}) define a conformal immersion
of surface into $\C^4$ with induced metric of the form
$$ds^2=\psi_1\psi_2\varphi_1\varphi_2dzd\bar{z}.$$

On the other hand, the formulae (\ref{e1}) may be rewritten as follows
\begin{eqnarray}
d(X^1+iX^2)&=&i\psi_1\psi_2d\bar{z},\nonumber\\
d(X^1-iX^2)&=&-i\varphi_1\varphi_2dz,\nonumber\\
d(X^4+iX^3)&=&i\psi_1\varphi_2d\bar{z},\nonumber\\
d(X^4-iX^3)&=&i\varphi_1\psi_2dz\nonumber
\end{eqnarray}
or
\begin{equation}\label{e3}
d(X^4\sigma_0+X^1\sigma_1+X^2\sigma_2+X^3\sigma_3)=i\begin{pmatrix}
\varphi_1\psi_2dz & \psi_1\psi_2\bar{z}\\
\varphi_1\varphi_2dz & \psi_1\varphi_2d\bar{z}
\end{pmatrix},
\end{equation}
where $\sigma_0=\begin{pmatrix} 1 & 0 \\ 0 & 1 \end{pmatrix},\;\sigma_1=
\begin{pmatrix} 0 & 1 \\ -1 & 0 \end{pmatrix},\;\sigma_2=\begin{pmatrix}
0 & i \\ i & 0 \end{pmatrix},\;\sigma_3=\begin{pmatrix} -i & 0 \\ 0 & i
\end{pmatrix}$ are the matrix representations of units of a quaternion
algebra $\cl_{0,2}=\BH$. Since the coordinates $X^i$ are complex, then it is
easy to see that the left part of expression (\ref{e3}) be a biquaternion
$\C_2=\C\otimes\cl_{0,2}$. Moreover, there is an identity $\C_2=\cl_{3,0}$.
Indeed, a general element of algebra $\cl_{3,0}$ defined by the following
expression
\begin{equation}\label{e4}
\cA=a^0\e_0+\sum^3_{i=1}a^i\e_i+\sum^3_{i=1}\sum^3_{j=1}a^{ij}\e_{ij}+
a^{123}\e_{123}.
\end{equation}
It is obvious that a volume element $\omega=\e_{123}$ commutes with all
basis elements of algebra $\cl_{3,0}$. Therefore, a center of $\cl_{3,0}$
consists of the unit $\e_0$ and element $\omega$. Moreover, $\omega^2=-1$,
it is allows to identify the volume element with imaginary unit $i$.
Further, using the obvious identities
\begin{eqnarray}
&&\omega\e_1=\e_1\omega=\e_{23},\nonumber\\
&&\omega\e_2=\e_2\omega=\e_{31},\nonumber\\
&&\omega\e_3=\e_3\omega=\e_{12}\nonumber
\end{eqnarray}
the element (\ref{e4}) may be rewritten in the form
\begin{multline}
\cA=a^0\e_0+a^1\e_1+a^2\e_2+a^3\e_3+a^{12}\e_{12}+a^{31}\e_{31}+a^{23}\e_{23}+
a^{123}\e_{123}=\\
=(a^0+\omega a^{123})\e_0+(a^1+\omega a^{23})\e_1+(a^2+\omega a^{31})\e_2+
(a^3+\omega a^{12})\e_3.
\end{multline}
Recalling that $i\equiv\omega=\e_{123}$ and suppose $\e_3=\e_2\e_1$ we
obtain $\cl_{3,0}=\C_2$, where $\C_2$ is an algebra of complex quaternions,
the general element of which has a form
$$\cA=c^0\e_0+c^1\e_1+c^2\e_2+c^3\e_{21},$$
where $c^i\in\C$. Using the identity $\cl_{3,0}=\C_2$ and denoting
$X^4=X^0$ we can rewrite the left part of expression (\ref{e3}) as follows
\begin{multline}\label{e5}
(\re X^0+\omega\im X^0)\sigma_0+(\re X^1+\omega\im X^1)\sigma_1+\\
+(\re X^2+\omega\im X^2)\sigma_2+(\re X^3+\omega\im X^3)\sigma_3=\\
=\re X^0\sigma_0+\re X^1\sigma_1+\re X^2\sigma_2+\re X^3\sigma_3+\\
+\im X^3\sigma_{12}+\im X^2\sigma_{31}+\im X^1\sigma_{23}+\im X^0\sigma_{123}.
\end{multline}

Let us consider a space-time algebra $\cl_{1,3}$. A general element of
$\cl_{1,3}$ has a form
$$
\cA=a^0+\sum^3_{i=0}a^i\Gamma_i+\sum^3_{i=0}\sum^3_{j=0}a^{ij}\Gamma_{ij}+
\sum^3_{i=0}\sum^3_{j=0}\sum^3_{k=0}a^{ijk}\Gamma_{ijk}+a^{0123}\Gamma_{0123},
$$
where
$$\Gamma_0=\begin{pmatrix}
1 & 0 & 0 & 0\\
0 & 1 & 0 & 0\\
0 & 0 &-1 & 0\\
0 & 0 & 0 &-1
\end{pmatrix},\quad\Gamma_1=\begin{pmatrix}
0 & 0 & 0 & 1\\
0 & 0 & 1 & 0\\
0 &-1 & 0 & 0\\
-1& 0 & 0 & 0
\end{pmatrix},$$
\begin{equation}\label{e5a}
\Gamma_2=\begin{pmatrix}
0 & 0 & 0 &-i\\
0 & 0 & i & 0\\
0 & i & 0 & 0\\
-i& 0 & 0 & 0
\end{pmatrix},\quad\Gamma_3=\begin{pmatrix}
0 & 0 & 1 & 0\\
0 & 0 & 0 &-1\\
-1& 0 & 0 & 0\\
0 & 1 & 0 & 0
\end{pmatrix}.
\end{equation}
It is easy to see that a set $\cl^+_{1,3}$ of all even elements of space-time
algebra (subalgebra of $\cl_{1,3}$) is isomorphic to the biquaternion:
$\cl^+_{1,3}\cong\cl_{3,0}$. Further, Dirac algebra $\C_4$ be a 
complexification of space-time algebra: $\C_4=\C\otimes\cl_{1,3}$. On the
other hand, the volume element $\omega=\e_{01234}\in\cl_{4,1}$ is belong to
a center of $\cl_{4,1}$ and $\omega^2=-1$, therefore we have an identity
$\C_4=\cl_{4,1}$.

Consider now an important  notion of a minimal left ideal of
Clifford algebra. Let $\cl_{p,q}(V,Q)$ be a Clifford algebra over a real
field $\R$, where $V$ is
a vector space endowed with nondegenerate quadratic form
$$Q=x^2_1+\ldots+x^2_p-\ldots-x^2_{p+q}.$$    
A minimal left (respectively right) ideal is a set of type $I_{p,q}=\cl_{p,q}
e_{pq}$ (resp. $e_{pq}\cl_{p,q}$), where $e_{pq}$ is a primitive
idempotent, i.e., $e^2_{pq}=e_{pq}$ and $e_{pq}$ cannot be represented as a
sum of two orthogonal idempotents, i.e., $e_{pq}\neq f_{pq}+g_{pq}$, where
$f_{pq}g_{pq}=g_{pq}f_{pq}=0,\;f^2_{pq}=f_{pq},\,g^2_{pq}=g_{pq}$.
\begin{theorem}[{\rm Lounesto \cite{17}}] A minimal left ideal of $\cl_{p,q}$ 
is of the
type $I_{p,q}=\cl_{p,q}e_{pq}$, where $e_{pq}=\frac{1}{2}(1+\e_{\alpha_1})
\ldots\frac{1}{2}(1+\e_{\alpha_k})$ is a primitive idempotent of $\cl_{p,q}$
and $\e_{\alpha_1},\ldots,\e_{\alpha_k}$ are commuting elements of the
canonical basis of $\cl_{p,q}$ such that $(\e_{\alpha_i})^2=1,\,(i=1,2,\ldots,
k)$ that generate a group of order $2^k,\;k=q-r_{q-p}$ and $r_i$ are the
Radon-Hurwitz numbers, defined by the recurrence formula $r_{i+8}=r_i+4$
and
\begin{center}
\begin{tabular}{lcccccccc}
$i$ & 0 & 1 & 2 & 3 & 4 & 5 & 6 & 7 \\ \hline
$r_i$ & 0 & 1 & 2 & 2 & 3 & 3 & 3 & 3
\end{tabular}.
\end{center}
\end{theorem} 
From adduced above theorem immediately follows that minimal left ideals
of space-time algebra $\cl_{1,3}$ and Dirac algebra $\cl_{4,1}$ have 
respectively a form
\begin{eqnarray}
I_{1,3}&=&\cl_{1,3}e_{13}=\cl_{1,3}\frac{1}{2}(1+\Gamma_0),\label{e5'}\\
I_{4,1}&=&\cl_{4,1}e_{41}=\cl_{4,1}\frac{1}{2}(1+\Gamma_0)\frac{1}{2}
(1+i\Gamma_{12}).\label{e5''}
\end{eqnarray}
Moreover, for the minimal left ideal of Dirac algebra using the isomorphisms
$\cl_{4,1}=\C_4=\C\otimes\cl_{1,3}\cong\M_2(\C_2),\;\cl^+_{4,1}\cong\cl_{1,3}
\cong\M_2(\BH)$ and also an identity $\cl_{1,3}e_{13}=\cl^+_{1,3}e_{13}$ 
\cite{14,18} we have the following expression \cite{19}
\begin{multline}\label{e6}
I_{4,1}=\cl_{4,1}e_{41}=(\C\otimes\cl_{1,3})e_{41}\cong\cl^+_{4,1}e_{41}\cong
\cl_{1,3}e_{41}=\\
=\cl_{1,3}e_{13}\frac{1}{2}(1+i\Gamma_{12})=\cl^+_{1,3}e_{13}\frac{1}{2}
(1+i\Gamma_{12}).
\end{multline}
Further, let $\Phi\in\cl_{4,1}\cong\M_4(\C)$ be a Dirac spinor and 
$\phi\in\cl^+_{1,3}\cong\cl_{3,0}=\C_2$ be a Dirac-Hestenes spinor. Then from
(\ref{e6}) immediately follows a relation between spinors $\Phi$ and $\phi$:
\begin{equation}\label{e7}
\Phi=\phi\frac{1}{2}(1+\Gamma_0)\frac{1}{2}(1+i\Gamma_{12}).
\end{equation}
Since $\phi\in\cl^+_{1,3}\cong\cl_{3,0}$, then the Dirac-Hestenes spinor can be
represented by a biquaternion number
\begin{multline}\label{e8}
\phi=a^0+a^{01}\Gamma_{01}+a^{02}\Gamma_{02}+a^{03}\Gamma_{03}+\\
+a^{12}\Gamma_{12}+a^{13}\Gamma_{13}+a^{23}\Gamma_{23}+a^{0123}\Gamma_{0123}.
\end{multline}
Or in the matrix representation
\begin{equation}\label{e9}
\phi=\begin{pmatrix}
\phi_1 & -\phi^\ast_2 & \phi_3 & \phi^\ast_4\\
\phi_2 & \phi^\ast_1 & \phi_4 & -\phi^\ast_3\\
\phi_3 & \phi^\ast_4 & \phi_1 & -\phi^\ast_2\\
\phi_4 & -\phi^\ast_3 & \phi_2 & \phi^\ast_1
\end{pmatrix},
\end{equation}
where
\begin{eqnarray}
\phi_1&=&a^0-ia^{12},\nonumber\\
\phi_2&=&a^{13}-ia^{23},\nonumber\\
\phi_3&=&a^{03}-ia^{0123},\nonumber\\
\phi_4&=&a^{01}+ia^{02}.\nonumber
\end{eqnarray}
According to (\ref{e5''}) and (\ref{e7}), (\ref{e9}) in the matrix
representation elements of minimal left ideal of Dirac algebra have a form
\begin{equation}\label{e9'}
\Phi=\begin{pmatrix}
\phi_1 & 0 & 0 & 0\\
\phi_2 & 0 & 0 & 0\\
\phi_3 & 0 & 0 & 0\\
\phi_4 & 0 & 0 & 0
\end{pmatrix}.
\end{equation} 
Thus, the elements of this minimal left ideal contain four complex, or
eight real parameters, which are just sufficient to define a Dirac spinor.

Let us return to a spinor type form of Weierstrass representation (\ref{e3})
and (\ref{e5}). It is easy to see that by force of $\cl_{3,0}\cong\cl^+_{1,3}$
and $\cl^{++}_{4,1}\cong\cl^+_{1,3}\cong\cl_{3,0}$ the right part of expression
(\ref{e5}) is isomorphic to the following biquaternion
\begin{multline}
\phi=\re X^0I+\re X^1\Gamma_{01}+\re X^2\Gamma_{02}+\re X^3\Gamma_{03}+\\
+\im X^3\Gamma_{12}+\im X^2\Gamma_{31}+\im X^1\Gamma_{23}+\im X^0\Gamma_{0123},
\end{multline} 
which may be rewritten in the matrix form (\ref{e9}) if suppose
\begin{eqnarray}
\phi_1&=&\re X^0-i\im X^3,\nonumber\\
\phi_2&=&\im X^2-i\im X^1,\nonumber\\
\phi_3&=&\re X^3-i\im X^0,\nonumber\\
\phi_4&=&\re X^1+i\re X^2,\label{e10}
\end{eqnarray}
So, we establish a relation between Weierstrass coordinates for surfaces
conformally immersed into $\C^4$ and Dirac-Hestenes spinors. Further,
in accordance with (\ref{e7}) it is easy to establish a relation with the
Dirac spinor $\Phi\in\M_4(\C)e_{41}$ treated as minimal left ideal of
$\cl_{4,1}=\C_4\cong\M_4(\C)$. Therefore, a Dirac field
$\Phi=(\phi_1,\,\phi_2,\,\phi_3,\,\phi_4)^T$ (which as known described an
electron in physics) may be expressed by means of relations (\ref{e10})
via the generalized Weierstrass formulae (\ref{e1}). In some sense it is allows
to consider the electron as a surface conformally immersed into 
4-dimensional complex space $\C^4$.
\section{Charge conjugation and antiautomorphism $\cA\longrightarrow
\widetilde{\cA^\star}$}
In Clifford algebra $\cl_{p,q}$ there exist four fundamental automorphisms
\cite{20,21}:

1) An automorphism $\cA\rightarrow\cA$.\\
This automorphism, obviously, be an identical automorphism of algebra 
$\cl_{p,q}$, $\cA$ is an arbitrary element of $\cl_{p,q}$.

2) An automorphism $\cA\rightarrow\cA^\star$.\\
In more details, for arbitrary element $\cA\in\cl_{p,q}$ there exist a
decomposition
$$\cA=\cA^\p+\cA^{\p\p},$$
where $\cA^\p$ is an element consisting of homogeneous odd elements, and
$\cA^{\p\p}$ is an element consisting of homogeneous even elements,
respectively. Then the automorphism $\cA\rightarrow\cA^\star$ is that
element $\cA^{\p\p}$ is not changed, and element $\cA^\p$ is changed the
sign:
$$\cA^\star=-\cA^\p+\cA^{\p\p}.$$
If $\cA$ is a homogeneous element, then
\begin{equation}\label{e11}
\cA^\star=(-1)^k\cA,
\end{equation}
where $k$ is a degree of element.

3) An antiautomorphism $\cA\rightarrow\widetilde{\cA}$.\\
The antiautomorphism $\cA\rightarrow\widetilde{\cA}$ be a reversion
of the element $\cA$, that is the substitution of the each basis element
$\e_{i_1i_2\ldots i_k}\in\cA$ by the element $\e_{i_ki_{k-1}\ldots i_1}$:
$$\e_{i_ki_{k-1}\ldots i_1}=(-1)^{\frac{k(k-1)}{2}}\e_{i_1i_2\ldots i_k}.$$
Therefore, for any $\cA\in\cl_{p,q}$ we have
\begin{equation}\label{e12}
\widetilde{\cA}=(-1)^{\frac{k(k-1)}{2}}\cA.
\end{equation}

4) An antiautomorphism $\cA\rightarrow\widetilde{\cA^\star}$.\\
This antiautomorphism be a composition of the antiautomorphism
$\cA\rightarrow\widetilde{\cA}$ with the automorphism $\cA\rightarrow
\cA^\star$. In the case of homogeneous element from formulae (\ref{e11})
and (\ref{e12}) follows
\begin{equation}\label{e13}
\widetilde{\cA^\star}=(-1)^{\frac{k(k+1)}{2}}\cA.
\end{equation}
It is obvious that $\;\;\widetilde{\!\!\widetilde{\cA}}=\cA,
\;(\cA^\star)^\star=\cA,\;
\widetilde{\left(\widetilde{\cA^\star}\right)^\star}=\cA.$

The antiautomorhism $\cA\rightarrow\widetilde{\cA^\star}$ is present
for us a most interest, since this antiautomorphism is closely related
with a charge conjugation in theory of electron \cite{21}. 
It is known \cite{21} that in the
matrix representation the antiautomorpphism $\cA\rightarrow\widetilde{
\cA^\star}$ defined by the following expression
\begin{equation}\label{e14}
\widetilde{A^\star}=(CE^T)A^T(CE^T)^{-1},
\end{equation}
where $E$ is a matrix of volume element of $\cl_{p,q}$, $C$ is a matrix
of antiautomorphism $\cA\rightarrow\widetilde{\cA}$. Let us now find
a matrix $CE^T$ of antiautomorphism $\cA\rightarrow\widetilde{\cA^\star}$
for the space-time algebra $\cl_{1,3}$. In the spinbasis (\ref{e5a}) the
matrix $E$ has a form
$$E=\Gamma_{0123}=\begin{pmatrix}
0 & 0 & -i & 0\\
0 & 0 & 0 &-i\\
-i& 0 & 0 & 0\\
0 & -i& 0 & 0
\end{pmatrix}.$$
Further, under action of antiautomorphism $\cA\rightarrow\widetilde{\cA}$
the units of $\cl_{p,q}$ are transfered into themselves: $\e_i\rightarrow
\e_i\;(i=1,\ldots,n=p+q)$. Therefore, in the case of spinbasis (\ref{e5a})
we have also: $\Gamma_i\rightarrow\Gamma_i$. On the other hand, in the
matrix representation for $\cA\rightarrow\widetilde{\cA}$ we have \cite{21}
\begin{equation}\label{e15}
A\longrightarrow CA^TC^{-1}.
\end{equation}
Transposition of the matrices (\ref{e5a}) gives
$$\Gamma^T_0=\Gamma_0,\quad\Gamma^T_1=-\Gamma_1,\quad\Gamma^T_2=\Gamma_2,\quad
\Gamma^T_3=-\Gamma_3.$$
Take into account the last relations we obtain from (\ref{e15})
\begin{eqnarray}
\Gamma_0\longrightarrow\Gamma_0=C\Gamma_0C^{-1},&&\Gamma_1
\longrightarrow\Gamma_1=-C\Gamma_1C^{-1},\nonumber\\
\Gamma_2\longrightarrow\Gamma_2=C\Gamma_2C^{-1},&&\Gamma_3\longrightarrow
\Gamma_3=-C\Gamma_3C^{-1},
\end{eqnarray}
or
\begin{eqnarray}
\Gamma_0C=C\Gamma_0,&&\Gamma_1C=-C\Gamma_1,\nonumber\\
\Gamma_2C=C\Gamma_2,&&\Gamma_3C=-C\Gamma_3.\label{e16}
\end{eqnarray}
It is easy to see that a matrix $C=\Gamma_{13}$ is satisfy to conditions
(\ref{e16}) and therefore be a matrix of antiautomorphism $\cA
\rightarrow\widetilde{\cA}$ for the spinbasis (\ref{e5}). Hence it
immediately follows that for a matrix of antiautomorphism $\cA
\rightarrow\widetilde{\cA^\star}$ we have
$$CE^T=\begin{pmatrix}
0 & 0 & 0 & -i\\
0 & 0 & i & 0\\
0 &-i & 0 & 0\\
i & 0 & 0 & 0
\end{pmatrix}.$$
Further, using the matrix representation (\ref{e14}) we find that an action
of antiautomorphism $\cA\rightarrow\widetilde{\cA^\star}$ on 
Dirac-Hestenes spinor field is expressed as follows
\begin{equation}\label{e17}
\widetilde{\phi^\star}=\begin{pmatrix}
\phi^\ast_1 & \phi^\ast_2 & -\phi^\ast_3 & -\phi^\ast_4 \\
-\phi_2 & \phi_1 & -\phi_4 & \phi_3 \\
-\phi^\ast_3 & -\phi^\ast_4 & \phi^\ast_1 & \phi^\ast_2 \\
-\phi_4 & \phi_3 & -\phi_2 & \phi_1
\end{pmatrix}.
\end{equation}
Whence, in accordance with relations (\ref{e7}) for a charge conjugated
Dirac spinor we obtain
\begin{multline}\label{e18}
\widetilde{\Phi^\star}=\widetilde{(\phi e_{41})^\star}=\widetilde{e_{41}^\star}
\widetilde{\phi^\star}=\\
\begin{pmatrix}
0 & 0 & 0 & 0\\
0 & 0 & 0 & 0\\
0 & 0 & 0 & 0\\
0 & 0 & 0 & 1
\end{pmatrix}\begin{pmatrix}
\phi^\ast_1 & \phi^\ast_2 & -\phi^\ast_3 & -\phi^\ast_4 \\
-\phi_2 & \phi_1 & -\phi_4 & \phi_3 \\
-\phi^\ast_3 & -\phi^\ast_4 & \phi^\ast_1 & \phi^\ast_2 \\
-\phi_4 & \phi_3 & -\phi_2 & \phi_1 
\end{pmatrix}=\begin{pmatrix}
0 & 0 & 0 & 0\\
0 & 0 & 0 & 0\\
0 & 0 & 0 & 0\\
-\phi_4 & \phi_3 & -\phi_2 & \phi_1
\end{pmatrix}.
\end{multline}
Therefore, {\it Dirac spinor (\ref{e9'}) treated as a minimal left ideal
of algebra $\cl_{4,1}=\C_4\cong\M_4(\C)$ under action of antiautomorphism
$\cA\rightarrow\widetilde{\cA^\star}$ is transfered into a minimal
right ideal of type (\ref{e18}).} 
Therefore, the field $\widetilde{\Phi^\star}=(-\phi_4,\,\phi_3,\,-\phi_2,\,
\phi_1)$ is expressed by means of (\ref{e10}) via the coordinates of
surface conformally immersed into $\C^4$. Thus, we have the following
\begin{theorem}
Under action of antiautomorphism $\cA\rightarrow\widetilde{\cA^\star}$
of space-time algebra $\cl_{1,3}$ the Dirac spinor
\[
\Phi=\phi e_{41}=\begin{pmatrix}
\phi_1 & 0 & 0 & 0\\
\phi_2 & 0 & 0 & 0\\
\phi_3 & 0 & 0 & 0\\
\phi_4 & 0 & 0 & 0
\end{pmatrix}
\]
which be a minimal left ideal of Dirac algebra $\cl_{4,1}=\C_4\cong\M_4(\C)$,
is transfered into a minimal right ideal (a charge conjugated spinor)
\[
\widetilde{\Phi^\star}=\widetilde{e^\star_{41}}\widetilde{\phi^\star}=
\begin{pmatrix}
0 & 0 & 0 & 0\\
0 & 0 & 0 & 0\\
0 & 0 & 0 & 0\\
-\phi_4 & \phi_3 & -\phi_2 & \phi_1
\end{pmatrix}
\]
of the same algebra. Here $\phi$ is a Dirac-Hestenes spinor,
$e_{41}=\frac{1}{2}(1+\Gamma_0)\frac{1}{2}(1+i\Gamma_{12})$ is a primitive
idempotent of $\C_4$, and
\begin{eqnarray}
\phi_1&=&\re X^0-i\im X^3,\nonumber\\
\phi_2&=&\im X^2-i\im X^1,\nonumber\\
\phi_3&=&\re X^3-i\im X^0,\nonumber\\
\phi_4&=&\re X^1+i\re X^2,\nonumber
\end{eqnarray}
where
\begin{eqnarray}
X^1&=&\frac{i}{2}\int_\Gamma(\psi_1\psi_2d\bar{z}-\phi_1\phi_2dz),\nonumber\\
X^2&=&\frac{1}{2}\int_\Gamma(\psi_1\psi_2d\bar{z}+\phi_1\phi_2dz),\nonumber\\
X^3&=&\frac{1}{2}\int_\Gamma(\psi_1\phi_2d\bar{z}-\phi_1\psi_2dz),\nonumber\\
X^0&=&\frac{i}{2}\int_\Gamma(\psi_1\phi_2d\bar{z}+\phi_1\psi_2dz)\nonumber
\end{eqnarray}
are generalized Weierstrass formulae for conformal immersion of surface
into $\C^4$.
\end{theorem}
\section{Integrable deformations and time evolution of Dirac field}
The system (\ref{e2}) is known in soliton theory as a
Davey-Stewartson II (DSII) linear problem \cite{22,23}. 
Integrable deformations of surface conformally
immersed into $\C^4$ are defined by an infinite hierarchy of nonlinear
differential equations associated with the system (\ref{e2}). This 
hierarchy is appear as compatibility conditions of (\ref{e2}) with systems
of the following form
\begin{equation}\label{e20}
\begin{array}{ccc}
\psi_{\alpha t_n}&=&A_n\psi_\alpha + B_n\varphi_\alpha,\\
\varphi_{\alpha t_n}&=&C_n\psi_\alpha+D_n\varphi_\alpha,
\end{array}\quad\alpha=1,2,
\end{equation}
where $t_n$ are new deformation variables and $A_n,\,B_n,\,C_n,\,D_n$ are
differential operators of $n$-th order. For example, in the case $n=3$
there is the system
\begin{eqnarray}
p_{t_3}&=&p_{zzz}+p_{\bar{z}\bar{z}\bar{z}}+3p_z\partial^{-1}_{\bar{z}}(pq)_z+
3p_{\bar{z}}\partial^{-1}_z(qp)_{\bar{z}}+3p\partial^{-1}_{\bar{z}}(qp_z)_z+
3p\partial^{-1}_z(qp_{\bar{z}})_{\bar{z}},\nonumber\\
q_{t_3}&=&q_{zzz}+q_{\bar{z}\bar{z}\bar{z}}+3q_z\partial^{-1}_{\bar{z}}(pq)_z+
3q_{\bar{z}}\partial^{-1}_z(pq)_{\bar{z}}+3q\partial^{-1}_{\bar{z}}(pq_z)_z+
3q\partial^{-1}_z(pq_{\bar{z}})_{\bar{z}}\nonumber\\ \label{e21}
\end{eqnarray}
and
\begin{eqnarray}
A_3&=&\partial^3_{\bar{z}}+3[\partial^{-1}_z(pq)_{\bar{z}}]\partial_{\bar{z}}+
3[\partial^{-1}_z(qp_{\bar{z}})_{\bar{z}}],\nonumber\\
B_3&=&-p\partial^2_z+p_z\partial_z-p_{zz}-3p[\partial^{-1}_{\bar{z}}(pq)_z],
\nonumber\\
C_3&=&-q\partial^2_{\bar{z}}+q_{\bar{z}}\partial_{\bar{z}}-q_{\bar{z}\bar{z}}-
3q[\partial^{-1}_{\bar{z}}(pq)_{\bar{z}}],\nonumber\\
D_3&=&\partial^3_z+3[\partial^{-1}_{\bar{z}}(pq)_z]\partial_z+3[\partial^{-1}_
z(pq)_{\bar{z}}],
\end{eqnarray}

It is obvious that deformation of $\psi_\alpha,\,\varphi_\alpha$ by means
of (\ref{e20}) induces deformations of coordinates $X^i(z,\bar{z},t_n)$ of
surface in $\C^4$. Moreover, according to (\ref{e10}) DSII-deformation
generates a deformation (time evolution) of Dirac-Hestenes spinor field and
in accordance with (\ref{e9'}) a time evolution of Dirac field
$\Phi=(\phi_1,\,\phi_2,\,\phi_3,\,\phi_4)^T$. The all chain of deformations
may be represented by the following scheme
\[
p,\,q\longrightarrow\psi_\alpha,\,\varphi_\alpha\longrightarrow X^i
\longrightarrow\phi_i\longrightarrow\Phi.
\]
\begin{prop} A time evolution of the Dirac fields $\Phi=(\phi_1,\,\phi_2,
\,\phi_3,\,\phi_4)^T$ and $\widetilde{\Phi^\star}=(-\phi_4,\,\phi_3,\,-\phi_2
,\,\phi_1)$, where the components $\phi_i$ are expressed via coordinates
of conformally immersed surface into complex space $\C^4$, is defined by
the DSII-hierarchy.
\end{prop}

In particular case $p=q$ and $p=\bar{p}$ the equations (\ref{e21}),
 reduce to a modified Veselov-Novikov equation \cite{24}
\[
p_t=p_{zzz}+p_{\bar{z}\bar{z}\bar{z}}+3p_z\partial^{-1}_{\bar{z}}(|p|^2_z)+
3p_{\bar{z}}\partial^{-1}_z(|p|^2_{\bar{z}})+3p\partial^{-1}_{\bar{z}}
(\bar{p}p_z)_z+3p\partial^{-1}_z(\bar{z}p_{\bar{z}})_{\bar{z}}.
\]
Therefore, {\it integrable deformations of surfaces in $\C^4$ with $p=q$,
$p=\bar{p}$ are defined by the mVN-hierarchy}.

In other particular case $q=1$ and $p$ is a real-valued function the
equations (\ref{e21}), reduce to Veselov-Novikov equation 
\cite{25,26}
\[
p_t=p_{zzz}+p_{\bar{z}\bar{z}\bar{z}}+3[p\partial^{-1}_{\bar{z}}(p_z)]_z+
3[p\partial^{-1}_z(p_{\bar{z}})]_{\bar{z}}.
\]
So, {\it integrable deformations of surfaces in $\C^4$ with $q=1$ are
generated by the Veselov-Novikov hierarchy}.

\begin{flushleft}
Vadim V. Varlamov\\
Siberia State University of Industry\\
Kirova 42\\
Novokuznetsk 654007\\
Russia\\
E-mail: root@varlamov.kemerovo.su
\end{flushleft}
\end{document}